\documentclass[11pt,reqno]{amsart}
\usepackage{amssymb,euscript,epsfig}
\newcommand{\e}{\varepsilon}
\newcommand{\eps}{\varepsilon}
\newcommand{\om}{\omega}

\newcommand{\la}{\lambda}
\newcommand{\al}{\alpha}

\newcommand{\fy}{\varphi}
\newcommand{\pd}{\partial}
\newcommand{\p}{\partial}

\newcommand{\ti}{\widetilde}
\newcommand{\R}{\mathbb{R}}
\newcommand{\C}{\mathbb{C}}

\newcommand{\Z}{\mathbb{Z}}
\newcommand{\ZZ}{\mathcal{Z}}
\newcommand{\F}{\mathcal{F}}

\newcommand{\NN}{\mathcal{N}}
\renewcommand{\Re}{\mathop{\mathrm{Re}}}
\renewcommand{\Im}{\mathop{\mathrm{Im}}}

\renewcommand{\bar}{\overline}
\renewcommand{\hat}{\widehat}
\newcommand{\FS}{{L^\infty\cap L^2}}
\newcommand{\WL}[3]{{L^{#1}_{#2,#3}}}

\numberwithin{equation}{section}

\newtheorem{theorem}{Theorem}[section]

\newtheorem{lemma}[theorem]{Lemma}

\theoremstyle{remark}
\newtheorem{remark}[theorem]{Remark}

\newcommand{\norm}[1]{\left \| #1 \right \|}
\newcommand{\ran}{\rangle}
\newcommand{\lan}{\langle}

\newcommand{\lec}{{\ \lesssim \ }}
\newcommand{\gec}{{\ \gtrsim \ }}

\newcommand{\sgec}{\gtrsim}

\newcommand{\EQ}[1]{\begin{equation} \begin{split} #1
 \end{split} \end{equation}}

\newcommand{\LR}[1]{{\lan #1 \ran}}
\newcommand{\Del}[1]{}
\newcommand{\de}{\delta}
\newcommand{\si}{\sigma}

\newcommand{\be}{\beta}
\newcommand{\ka}{\kappa}
\newcommand{\x}{\xi}
\newcommand{\y}{\eta}
\newcommand{\z}{\zeta}
\newcommand{\s}{\sigma}
\newcommand{\rh}{\rho}
\newcommand{\na}{\nabla}
\renewcommand{\th}{\theta}

\newcommand{\supp}{\operatorname{supp}}
\newcommand{\Tri}{{\it Tri}}
\newcommand{\Dif}{{\it Dif}}
\newcommand{\Asy}{{\it Asy}}
\newcommand{\intS}[1]{\int \int_\infty^t e^{i\Phi_+ s}(\chi_{X}^{#1}F) ds d\y}
\newcommand{\X}{\mathcal{X}^}

\newcommand{\pp}{^\perp_}
\newcommand{\BR}[1]{[#1]}

\begin{document}

\title[Global dispersive solutions for Gross-Pitaevskii in 2D \& 3D]{Global dispersive solutions for the Gross-Pitaevskii equation\\ in two and three dimensions}
\author{Stephen Gustafson,\quad Kenji Nakanishi,\quad Tai-Peng Tsai}

\begin{abstract}
We study asymptotic behaviour at time infinity of solutions close to the non-zero constant equilibrium for the Gross-Pitaevskii equation in two and three spatial dimensions. We construct a class of global solutions with prescribed dispersive asymptotic behavior, which is given in terms of the linearized evolution.  
\end{abstract}

\maketitle

\section{Introduction}
We consider the Gross-Pitaevskii equation for $\psi:\R^{1+d}\to\C$ 
\EQ{\label{GP}
  i \pd_t \psi = - \Delta \psi + (|\psi|^2-1)\psi,}
with the boundary condition
\EQ{ \label{BC}
 |\psi(t,x)| \to 1\quad \mbox{ as } \quad |x| \to \infty.}
This equation is a dynamical model for superfluids and Bose-Einstein condensates, and has been extensively studied, especially concerning traveling wave solutions of the form $\psi=\fy(x-ct)$, and dynamics of vortices (zeros of $\psi$). See [2--10, 12--16, 20--25, 27, 28] and references therein. 
However we know very little about long-time dynamics of general
solutions, for example about stability of vortices and traveling
waves, or even of the constant solution $\psi=1$. Heuristically, the main difficulty is that small perturbations can continue to interact with the non-zero background $1$, and so do not easily disperse and decay. 

Thus we started in \cite{vac} an investigation of large-time behavior of solutions $\psi=1+$``small" as a first step toward understanding dispersive processes in this equation. The perturbation $u=\psi-1$ from the equilibrium satisfies the equation 
\EQ{ \label{eq u}
\begin{split}
&i \pd_t u + \Delta u - 2\Re u = F(u), \quad 
F(u) := u^2 + 2|u|^2 + |u|^2 u.
\end{split}}
The conserved energy and charge are written respectively 
\EQ{
 \int_{\R^d} |\na u|^2 + \frac{(|u|^2+2\Re u)^2}{2}\ dx,\quad \int_{\R^d} |u|^2 + 2\Re u\ dx.} 
For $d=2$ and $3$, unique global existence for the Cauchy problem has been proved first in \cite{BS} for $u\in H^1$, and later in \cite{Ge} for any finite energy solution. 

In order to investigate dispersive properties of small solutions $u$, it is natural to linearize the equation around $0$. The left hand side of \eqref{eq u} can be made complex linear by the following change of variable: 
\EQ{
 u \mapsto v := V^{-1}u := U^{-1}\Re u + i \Im u, \quad U:=\sqrt{-\Delta(2-\Delta)^{-1}}.}
Then the new function $v$ satisfies the equation 
\EQ{ \label{eq v}
 i \p_t v - H v = -i V^{-1} i F(Vv),\quad H:=\sqrt{-\Delta(2-\Delta)}.}
The linear evolution $v=e^{-iHt}v(0)$ is expected to approximate small perturbations from the equilibrium. 
We have proved in the previous paper \cite{vac} that this is indeed the case if $d=4$ for small solutions $v\in H^1$.

 In this paper we turn to the physical dimensions $d=2,3$, where the dispersion becomes much weaker. Actually there exist traveling wave solutions with arbitrarily small energy for $d=2$ \cite{BS2}, so it seems unlikely that the same result holds as in $d=4$. However, it is still possible that small solutions $u$ disperse in general if they are well localized in space at some time, since those traveling waves are spatially spread-out and do not belong to $L^2$ \cite{Gr3}. Our theorems \ref{thm:3d}, \ref{thm:2d} show that at least there exist plenty of global dispersive solutions even for $d=2,3$. 

Since \eqref{eq v} is rather complicated, it seems natural to compare it with a simpler nonlinear Schr\"odinger equation (NLS) with the same nonlinearity:
\EQ{ \label{NLS}
 i\p_t v + \Delta v = F(v).}
The scaling argument suggests that the quadratic nonlinear terms can be treated as a perturbation in $L^\infty_t(L^2_x)$ only if $d\ge 4$. Hence we are led to work in weighted spaces, which provides more decay in time. Even with the optimal decay rate of the free evolution, the quadratic terms exhibit in general the critical decay order $1/t$ in $L^2_x$ if $d=2$. For $d=3$, we can generally expect asymptotically free behavior for dispersive solutions. 

Coming back to our equation, the operator $H$ has a singularity at 0 frequency $\x=0$ similar to the wave equation, which is worse for time decay than the Schr\"odinger equation. In addition, \eqref{eq v} apparently contains a singularity due to $V^{-1}$ in the nonlinearity. 
Nevertheless, it turns out that our equation is better than the NLS \eqref{NLS}, and in fact our argument for $d=2$ does not apply to the latter, which appears simpler at first glance. This is because \eqref{eq v} has a special structure and thereby we can transform it to another equation with a derivative nonlinearity, in effect. We give the details of the transform in Section \ref{ss:normal}.  

Before stating our results, we review the known results on the NLS for comparison. Dispersive global solutions have been constructed for the quadratic NLS in $d=2$ only in the following two cases: 
\begin{enumerate}
\item $F(u)=\la_0 |u|u + \la_1 u^2 + \la_2 \bar{u}^2$, ($\la_j\in\C$), \cite{GO,HNST} \label{case1}
\item $F(u)=[\Re(\la u)]^2$, ($\la\in\C$), \cite{HN}
\end{enumerate}
In the first case, the asymptotic profile is modified from the free
evolution by a phase factor which depends only on $\la_0|u|u$, while
in the second case it is modified by the addition of a term with $L^2$
scaling concentration at $\x=0$. Each argument depends essentially on
the form of the modification, and so it seems difficult to combine
these results to cover our $F(u)$. For example when $F(u)=|u|^2$, it
is only known that there are no asymptotically free solutions in the
$L^2$ sense with the natural decay property \cite{Sh}.

For $d=3$, it is known \cite{HMN} that small initial data in certain weighted spaces lead to asymptotically free solutions in the above case \ref{case1}. The final data problem (or construction of the wave operators) is easier and arbitrary quadratic terms can be treated in the same way as in Section \ref{s:3d}. 

Now we state the main results in this paper. $H^s_p$ and $\dot H^s_p$ denote the inhomogeneous and homogeneous Sobolev spaces respectively (cf. \cite{BL}), 
and we omit the subscript when $p=2$. Denote by $\dot B^s_{p,q}$ the homogeneous Besov spaces.
In three dimensions, we have wave operators without size restriction. 
\begin{theorem} \label{thm:3d}
Let $d=3$ and $0<\eps$ be small ($\eps\le 3/68$ is sufficient). For any $T\ge 1$, we define a Banach space $X^\eps_T$ by the following norm
\EQ{
 \|u\|_{X^\eps_T} = \sup_{T\le S} S^{1/2-8\eps}\|u\|_{L^p_t(S,\infty;H^1_q)},}
where $1/p=10\eps$ and $1/q=1/3-\eps$. For any $\fy\in H^1$ satisfying $\|e^{-iHt}\fy\|_{X^\eps_1} < \infty$, 
there exists a unique global solution $\psi=1+u$ of \eqref{GP} satisfying 
\EQ{
 &C(\R;H^1)\ni V^{-1}u = e^{-iHt}\fy + v',\\  
 &\|v'(t)\|_{(L^\infty_t H^1_x \cap L^2_t H^1_6)(T,\infty)} \lec T^{-1/4-\eps},\quad \|v'\|_{X^\eps_T} \lec T^{-\eps/2}.}
The above condition on $\fy$ is satisfied if $\fy\in H^1 \cap H^1_{q/(q-1)}$. 
We have the same result in the critical case $\eps=0$ if $\|e^{-iHt}\fy\|_{X^0_1}$ is small enough.  
\end{theorem}
The threshold $\eps=0$ is related to the scaling property of the NLS with quadratic nonlinearity in $d=3$. 

In two dimensions, we can construct asymptotically free solutions for small final data. 
\begin{theorem} \label{thm:2d}
Let $d=2$ and $\fy\in H^1$. Assume that 
\EQ{ \label{asy norm}
 \LR{\x}^{-1/2}|\x|^{|k|}\p^k \F \fy(\x) \in \FS}
for all multi-indices $k\ge 0$ with $|k|=k_1+k_2\le 2$, 
and that $\|\fy\|_{\dot B^1_{1,1}}$ is sufficiently small. Then there exists a unique global solution $\psi=1+u$ for \eqref{GP} satisfying 
\EQ{
 &V^{-1}u = z^{0} -\nu + z' + z'' \in C(\R;\dot H^1\cap\dot H^\eps),\\ 
 &z^{0} = e^{-iHt}\fy,\quad  
 z' = i\int_\infty^t e^{-iH(t-s)}|Uz^{0}|^2 ds,\quad
  \nu = (2-\Delta)^{-1}U^{-1}|u|^2,\\ 
 &\|z'\|_{\dot H^1}
  +\|z''\|_{H^1} +\|\nu\|_{\dot H^2\cap \dot H^1} 
  \lec t^{-1+\eps},\quad  
 \|z'\|_{\dot H^\eps}+\|\nu\|_{\dot H^\eps} \lec t^{-\eps/2}
}
for any $\eps\in(0,1)$ and $t>0$, where the constants depend on $\eps$. 
\end{theorem} 
\begin{remark}
The correction term $\nu$ is coming from the normal form (see Sect. \ref{ss:normal}). It has a singularity at $\x=0$, which can be worse than $1/|\x|$, because we do not know whether our solution $u$ belongs to $L^2_x$. That is also the reason we describe it in terms of $u$, not $\fy$. 

The correction term $z'$ is essentially the same as in \cite{HN} for the NLS with $(\Re u)^2$, although we do not know whether it can be simplified as there, because of the singularity of our $H(\x)$ at $\x=0$. It is probably not in $L^2_x$ in general.

However, these correction terms have no essential effect in the nonlinearity, and so they can be regarded as error terms if one does not require $L^2$ asymptotics for $v$ or $u_2$. 
\end{remark}

In the next subsections, we explain our basic tools, namely the normal form and the $L^p$ decay estimate. 

\subsection{Normal form} \label{ss:normal} 
To eliminate the singularity at zero frequency $\x=0$, we introduced in \cite{vac} the following transformation of normal form type:
\EQ{
 w = u + P\frac{|u|^2}{2},}
where $P$ was a Fourier multiplier cutting-off the higher frequency $|\x|\gec 1$. The new function $w$ satisfies the following equation
\EQ{
 i\dot w = &-\Delta w + 2\Re w + G(u),\\
  G(u) = &(3-P)u_1^2 + Q u_2^2 + P\Delta|u|^2/2 + |u|^2 u_1,\\ 
  &+i[2Q(u_1 u_2) + \nabla P\cdot(u_2\nabla u_1 - u_1\nabla u_2) +Q(|u|^2u_2)],}
where $u=u_1+iu_2$ and $Q=Id-P$. It was crucial in \cite{vac} for
$d=4$ that $\Im G$ is essentially of derivative form. We also
exploited the fact that the quadratic part does not contain $u_2^2$ in the low frequency. 

Here we make a new observation that a special choice of $P$ related to the equation leads to even better and much simpler nonlinearity. Let $Q=U^2=-\Delta/(2-\Delta)$ and $P=1-U^2=2/(2-\Delta)$. Then we have $2Q=-P\Delta$ and 
\EQ{
 G(u) = 2u_1^2 +|u|^2 u_1 - 2i\na\cdot P(u_1\na u_2) + iQ(|u|^2u_2).}
Hence the equation for $z=V^{-1}w=U^{-1}(u_1+P|u|^2/2)+iu_2$ is given by
\EQ{ \label{eq nu}
 &z = z^0 + \int_\infty^t e^{-iH(t-s)}[N^2(u)+N^3(u)] ds,\\
 &u = Vz - P\frac{|u|^2}{2},\quad z^0=e^{-iHt}\fy,}
where we denote
\EQ{
 N^2(u):=-2iu_1^2-2PU^{-1}\na\cdot(u_1\na u_2),\quad
 N^3(u):=-i|u|^2u_1 + U(|u|^2u_2).}
The new nonlinearity is roughly of the form $(Uz)^2 + U(z^3)$. It is vital for our analysis in $d=2$ that the quadratic terms consist only of derivatives. 

We will solve the above equation \eqref{eq nu} for $(z,u)$ and for $t>T\gg 1$ by the fixed point argument. Then solving the equation \eqref{eq u} for $u$ from $t=T$ by the result in \cite{Ge} (or by \cite{BS} for $d=3$), and using local uniqueness of $(z,u)$ satisfying 
\EQ{ \label{eq nu T}
 &z = e^{-iH(t-T)}z(T) + \int_T^t e^{-iH(t-s)}[N^2(u)+N^3(u)] ds,\quad 
 u = Vz - P\frac{|u|^2}{2},}
we can deduce that our solution $\psi:=u+1$ satisfies the Gross-Pitaevskii equation \eqref{GP} and extends globally in time. 
 
\subsection{$L^p$ decay estimate} 
We recall the linear decay estimate proved in \cite{vac}. We call the pair of exponents $(p,q)$ admissible if $2\le p,q\le\infty$, $(p,q)\not=(2,\infty)$ and $2/p + d/q = d/2$. We denote by $q'=q/(q-1)$ the H\"older conjugate. 
\begin{lemma} \label{th:4-2}
Let $d\ge 2$. (i) Let $2\le q\le \infty$ and $\s=1/2-1/q$. Then we have
\EQ{ \label{decay est}
 \|e^{-itH}\fy\|_{\dot B^0_{q,2}} \lec t^{-d\si} \|\fy\|_{\dot B^0_{q',2}}.
}
(ii) Let $(p,q)$, $(p_1,q_1)$ and $(p_2,q_2)$ be admissible. Then we have
\EQ{ \label{Strichartz}
 &\|e^{-itH}\fy\|_{L^{p} \dot B^0_{q,2}} \le C(p)^2 \|\fy\|_{L^2},\\
 &\left\|\int_{-\infty}^t e^{-i(t-s)H}f(s)ds\right\|_{L^{p_1} \dot B^0_{q_1,2}}
 \le C(p_1)C(p_2)\|f\|_{L^{p_2'} \dot B^0_{q_2',2}},
}
where $C(p)$ is some positive continuous function of $p$, but diverges as $p\to 2$ when $d=2$. 
\end{lemma}
The above estimates are exactly the same as for the Schr\"odinger evolution $e^{it\Delta}$. We had in \cite{vac} some gain at $\x=0$ for $d=3$, but we ignore it in this paper. The second last statement in Theorem \ref{thm:3d} follows from the above estimate (i). 

For any $s\in\R$ and $T\in\R$, we denote the full set of Strichartz norms of $H^s$ solutions for $t>T$ by 
\EQ{
 \|u\|_{Stz^s_T} := \sup_{(p,q):admissible} C(p)^{-2}\|u\|_{L^p(T,\infty;H^s_q)}.} 
When $d=3$, this is just $L^\infty H^s\cap L^2 H^s_6$. When $d=2$, it is slightly bigger than $L^\infty H^s\cap L^2H^s_\infty$. We define the weighted Lebesgue space $\WL{s}{b}{T}$ by the following norm for any $0\le b\le 1$, $s\in\R$ and $T>0$: 
\EQ{
 \|u\|_{\WL{s}{b}{T}} := \sup_{T\le S} S^{s}\|u(t)\|_{L^{1/b}_t(S,2S)}.}
We denote the mixed norm by (where $B$ is a Banach space) 
\EQ{
 \|u\|_{\WL{s}{b}{T}(B)} := \norm{\|u(t)\|_{B_x}}_{\WL{s}{b}{T}}.}
The H\"older inequality implies that 
\EQ{
 \WL{s_1}{b_1}{T} \times \WL{s_2}{b_2}{T} \subset \WL{s_1+s_2}{b_1+b_2}{T}.}
We have also $\WL{s_1}{b_1}{T}\subset\WL{s_2}{b_2}{T'}$ iff
\EQ{
 b_1\le b_2,\quad s_1+b_1\ge s_2+b_2,\quad T\le T'.}
Moreover we have 
\EQ{
 t^{-s} \in \WL{s}{0}{T} \quad(T>0).}

The rest of this paper is organized as follows. In the Section 2, we deal with the three dimensional case, and the other sections are devoted to two dimensions. After explaining the main ideas in Section 3, we give the main bilinear estimate in Section 4, and then prove Theorem \ref{thm:2d} in Section 5. 

\section{Three dimensions} \label{s:3d}
In this section, we construct the wave operators in $d=3$. 
The nonlinear terms are estimated simply by the H\"older and Sobolev inequalities, and the wave operators are constructed for the equation in our normal form by the standard fixed point theorem using the linear decay estimate. $H^{1/2}$ regularity would be sufficient for the final state problem, but we do not pursue it in this paper. Once the solution $u$ is constructed in $C([T,\infty);H^1)$ for some large $T>1$, it is uniquely extended to a global one by the result in \cite{BS}. We will construct the asymptotically free solution by the fixed point theorem in the space
\EQ{ \label{space 3D} 
  (z,u) \in Stz_T^1 \cap X^\eps_T,}
for large $T\gg 1$. 

\subsection{The scaling critical case}
We start with the simpler critical case $\eps=0$. By using the $L^3$ decay estimate, we have for the quadratic term for $t>T$, 
\EQ{
 \int_\infty^t \|e^{iH(s-t)}N^2(u)\|_{H^1_3} ds
 &\lec \int_\infty^t |s-t|^{-1/2}s^{-1}\|s^{1/2}u(s)\|_{H^1_3}^2 ds\\
 &\lec \|u\|_{X^0_T}^2 t^{-1/2},}
and for the cubic term 
\EQ{
 \int_\infty^t \|e^{iH(s-t)}N^3(u)\|_{H^1_3} ds
 &\lec \int_\infty^t |s-t|^{-1/2} s^{-1} \|s^{1/2}u(s)\|_{H^1_3}^2 \|u(s)\|_{L^\infty} ds\\
 &\lec t^{-3/4}\|u\|_{X_T^0}^2 \|u\|_{\WL{0}{1/4}{T}(L^\infty)}
  \lec t^{-3/4}\|u\|_{X_T^0}^2 \|u\|_{Stz^1_T}.}
The decay in $Stz^1$ is derived for the quadratic terms by 
\EQ{
 \norm{\int^t_\infty e^{iH(s-t)}N^2(u) ds}_{Stz^1_T}
  &\lec \|N^2(u)\|_{\WL{0}{3/4}{T}(H^1_{3/2})}\\ 
  &\lec \|t^{-1}\|t^{1/2}u(t)\|_{H^1_3}^2\|_{\WL{0}{3/4}{T}}
  \lec T^{-1/4}\|u\|_{X^0_T}^2 ,}
and for the cubic terms by 
\EQ{ \label{cubic est}
 \|N^3(u)\|_{\WL{0}{3/4}{T}(H^1_{3/2})}
 &\lec \|t^{-1}\|t^{1/2}u(t)\|_{H^1_3}^2\|u(t)\|_{L^\infty}\|_{\WL{0}{3/4}{T}} \\ 
 &\lec T^{-1/2}\|u\|_{X^0_T}^2\|u\|_{Stz^1_T}.}
As for the normal form, we have 
\EQ{
 \|U^{-1}P(\fy\psi)\|_{H^1_p} \lec \|\fy\|_{L^3} \|\psi\|_{L^p},}
for any $p>3/2$. Thus we get the unique solution $(z,u)$ for \eqref{eq nu} by the standard fixed point argument in the space \eqref{space 3D}, provided that $\|e^{-iHt}\fy\|_{X_T^0}$ is sufficiently small and $T>1$. Then the solution is extended globally by the result in \cite{BS}, and local uniqueness of the solution $(z,u)$ in $C_t(H^1)$ for \eqref{eq nu T}, which follows easily from the Strichartz, Sobolev and H\"older inequalities. 

\subsection{Large data wave operators}
Next we consider the case $\eps>0$ without size restriction. We define exponents $q,q^2,q_2$ by 
\EQ{
 1/q=1/3-\eps,\ 1/q^2=2/3-2\eps,\ 1/q_2=1/3+2\eps.}
The decay estimate implies that for $t>T$, 
\EQ{
 \norm{\int_\infty^t e^{-iH(t-s)}N^2(u) ds}_{H^1_{q_2}} &\lec \int_\infty^t |t-s|^{-1/2+6\e} \|N^2(u(s))\|_{H^1_{q^2}} ds\\
 &\lec T^{-1/2+2\e}\|u\|_{\WL{1/2-8\eps}{10\eps}{T}(H^1_q)}^2.}
The Strichartz estimate implies that 
\EQ{
 \norm{\int_\infty^t e^{-iH(t-s)}N^2(u) ds}_{Stz^1_T}
 &\lec \|N^2(u)\|_{\WL{0}{3/4+3\eps}{T}(H^1_{{q^2}})}
 \lec T^{-1/4-\eps}\|u\|_{\WL{1/2-8\eps}{10\eps}{T}(H^1_q)}^2.}
The cubic term has additional $T^{1/4}$ decay due to the $L^4_tL^\infty_x$ bound by the same argument as in \eqref{cubic est}. Then we use the complex interpolation   
\EQ{
 [\WL{1/2-2\e}{0}{T}(H^1_{q_2}), \WL{1/4+\eps}{1/2}{T}(H^1_6)]_\th = \WL{1/2-\th/4-\e(2-3\th)}{\th/2}{T}(H^1_{q}),}
where $\th\in[0,1]$ should be chosen to satisfy
\EQ{
 (1-\th)/q_2 + \th/6 = 1/q,}
i.e., $\th=18\e/(1+12\e)\le 1$. The last inequality is because $\eps\le 1/6$. To embed the above space into $X^\eps_T$, we need
\EQ{
 \th/2\le 10\e,\quad \th/4-\e(2-3\th) \ge 2\e,}
i.e., $16\e - 12\e\th \le \th \le 20\e$, which is satisfied by the above $\th$ with strict inequalities. In fact we have
\EQ{
 \th/4 - \eps(2-3\th) = 2\eps + 5\eps/2}
Therefore we get $T^{-5\eps/2}$ as a small factor for the nonlinear term in the space \eqref{space 3D}. The rest of proof is the same as in the critical case. \qedsymbol 

\section{Main ideas in Two dimensions}
In the rest of the paper, we deal with the case $d=2$. In this section we describe the outline, and derive the key estimate in the next section, then finally prove the main theorem in the last section. 

\subsection{Iteration scheme} 
Let $u^0:=Vz^0=Ve^{-iHt}\fy$. The integral equation is decomposed as follows
\EQ{
 z-z^0 = &\int_\infty^t e^{-iH(t-s)}N^3(u) ds\\
 &+\int_\infty^t e^{-iH(t-s)}[N^2(u)-N^2(u^0)] ds\\
 &+\int_\infty^t e^{-iH(t-s)}N^2(u^0) ds =: \Tri(u)+\Dif(u)+\Asy(u^0)}
The first two terms are estimated by simple H\"older and Sobolev type inequalities, and the main task is to derive enough time decay for the last term $\Asy(u^0)$, which is explicitly given by the data $\fy$. In estimating $\Dif(u)$, we use $t^{-1}L^\infty_x$ decay of $u^0$, which forces us to assume smallness of the data $\fy$ (this is usual in the case of critical decay). We further decompose $\Asy(u^0)$ as follows: 
\EQ{ \label{def Asy}
 \Asy(u^0) = \Asy'(u^0) + z',\quad z'=\int_\infty^t ie^{-iH(t-s)}|Uz^0|^2 ds.}
$z'$ is the only part where the oscillation of $u^0$ is completely canceled at $\x=0$. 

\subsection{Bilinear decay estimate}
We will derive in the next section 
\EQ{ \label{Asy decay}
 |\p^k\ti\fy| \lec |\x|^{-|k|}\LR{\x}^{1/2} \implies \Asy(u^0) \in t^{-1}(\log t)^2 \dot H^1.}
Notice that we have by the simple $t^{-1}L^\infty$ decay and the H\"older that 
\EQ{
 \fy \in \dot H^1\cap B^0_{1,1} \implies \int_{2t}^t e^{-iH(t-s)}N^2(u^0) ds \in L^\infty(H^1),}
so in \eqref{Asy decay} we are gaining roughly $1/t$ decay by losing $\x$ at $\x=0$, which is acceptable for our nonlinearity. 

The main idea of the decay estimate is as follows. For simplicity, consider the Schr\"odinger evolution $H(\x)=|\x|^2$. Our quadratic terms are roughly of the form
\EQ{ \label{quad}
 \int_\infty^t \int e^{i\Phi s} |\x|\ti\fy(\x-\y)\ti\psi(\y) d\y ds}
in the Fourier space, where the phase function $\Phi$ is given by one of
\EQ{
 &\Phi_0:=H(\x)-H(\x-\y)+H(\y)=|\x|^2-|\x-\y|^2+|\y|^2,\\
 &\Phi_\pm:=H(\x)\mp(H(\x-\y)+H(\y))=|\x|^2 \mp(|\x-\y|^2+|\y|^2).}
$\Phi_0$ corresponds to $|u|^2$ and $\Phi_\pm$ to $u^2$ and
$\bar{u}^2$.  We can gain $1/t$ by integration by parts in $\y$, picking
up the divisor $1/|\na_\y\Phi|$, where
\EQ{
 &\na_\y\Phi_0 = -(\y-\x)+\y = \x,\\
 &\na_\y\Phi_\pm=\mp[(\y-\x)+\y]=\mp(2\y-\x).}
Hence the singularity of $1/|\na_\y\Phi_0|$ is canceled by $|\x|$ in \eqref{quad}. We need to integrate twice, since we want to have $t^{-1}$ after the integration in $s$. Then we get $|\x|^{-1}$ in the case of $\Phi_0$, but it is almost in $L^2_\x$ and so OK if we allow the loss of $\log t$. 

In the case of $\Phi_\pm$, $\na_\y\Phi$ depends on $\y$, which reflects the fact that $u^2$ and $\bar{u}^2$ are oscillatory. But now we can integrate in $s$, because at the stationary point $\y=\x/2$, the phases $\Phi_\pm$ do not vanish:
\EQ{
 \Phi_+ = 2\y(\x-\y)= -2(\y-\x/2)^2-|\x|^2/2.}
We are getting strong divisor $1/|\x|^2$, but it is still OK around the stationary point, where we have $|\x|\sim|\y|\sim|\x-\y|$ and the nonlinearity supplies $|\y||\y-\x|$ decay. 

Since our actual symbol $H(\x)$ is degenerate at $\x=0$, we get a stronger singularity, where the $|\y||\y-\x|$ gain plays a crucial role. In addition, we should carefully compare the unbalanced radial and angular components. The detail starts in the next section. 

\section{Bilinear space-time phase estimates} \label{s:Asy}
As seen above, we are going to have a non-stationary phase estimate for the bilinear expression with integration in space-time. Now we state the main estimate in a slightly more general setting. 
\begin{lemma}
Let $d=2$ and $\s\ge -1/2$. Assume that $F(\x,\y)$ satisfies 
\EQ{ \label{cond F}
 &|\p_\y^{k} F| \lec \frac{|\y| |\x-\y|}{\LR{\x}\LR{\x-\y}^{\s}\LR{\y}^\s m^{|k|}\LR{m}} f(\x-\y)g(\y),\quad 
  |F| \lec f'(\x-\y)g'(\y),}
for all $0\le |k|\le 2$ and some nonnegative functions $f,g\in\FS$ and $f',g'\in L^2$, where $m:=\min(|\x-\y|,|\y|)$. $(f'$ and $g'$ are not related to $f$ and $g.)$ Then we have
\EQ{ \label{main est}
 &\norm{|\x|^\mu\int_\infty^t\int e^{i\Phi s} F(\x,\y)d\y ds}_{L^2_\x}\\
 & \lec t^{-\th}(\log t)^{1+\th}(\|f\|_{\FS}\|g\|_{\FS}+\|f'\|_{L^2}\|g'\|_{L^2}),}
for $0<\th\le 1$, $t >2$,
\EQ{
 \Phi=
 \begin{cases}
  \Phi_0 = H(\x)+H(\y)-H(\y-\x),\\
  \Phi_\pm = H(\x) \mp (H(\y)+H(\y-\x)),
 \end{cases}}
and 
\EQ{
 \begin{cases}
  \th\le\mu\le\th+(1+\s)\frac{1+\th}{2} &(\Phi=\Phi_0),\\
  \frac{\th-1}{2}\le\mu \le \th+(1+\min(\s,2\s))\frac{1+\th}{2} &(\Phi=\Phi_\pm).
 \end{cases}}
\end{lemma}

The main part of proof is to derive precise lower bounds on the first derivative of the phase and compatible upper bounds for the higher derivatives. In doing that, we should carefully distinguish the radial and angular components, otherwise we would get too much singularity at $\x=0$. 

\subsection{Preliminaries}
For any vectors $\x,\y$, we denote
\EQ{
 \LR{\x}=\sqrt{1+|\x|^2},\quad [\x]=\sqrt{2+|\x|^2},\quad \hat{\x}=\frac{\x}{|\x|},\quad 
 \x_\y:=\x\cdot\hat{\y},\quad \x\pp\y=\x-\hat{\y}\x_\y.}
Then the phase function $H(\x)$ is written as 
\EQ{
 H(\x) = H(|\x|) = |\x|[\x].}
We will denote $H'(\x):=H'(|\x|)$, etc. First we need to see that the above lemma applies to $\Asy(u^0)$. Its Fourier transform is a linear combination of the form $MG$ with 
\EQ{
 &M(\x,\y):=\begin{cases}
  |\y||\y-\x|[\y]^{-1}[\y-\x]^{-1} &(\text{for }u_1^2),\\
  \hat{\x}\cdot|\y|(\y-\x)[\x]^{-1}[\y]^{-1} &(\text{for }PU^{-1}\na\cdot(u_1\na u_2)),
 \end{cases}\\ 
 &G(\x,\y):=\begin{cases}
  \F\fy(\x-\y)\bar{\F\fy}(-\y) &(\Phi=\Phi_0=H(\x)+H(\y)-H(\y-\x)),\\
  \F\fy(\x-\y)\F\fy(\y) &(\Phi=\Phi_+=H(\x)-H(\y)-H(\y-\x)),\\
  \bar{\F\fy}(-\x+\y)\bar{\F\fy}(-\y) &(\Phi=\Phi_-=H(\x)+H(\y)+H(\y-\x)).
 \end{cases}}
Then in all six cases, our assumption \eqref{asy norm} implies the first condition of \eqref{cond F} with $\s=-1/2$ and some nonnegative functions $f,g\in\FS$ determined by $\fy$, and the second one follows from the assumption $\fy\in H^1$. In addition, we observe that $\Asy'(u^0)$ does not contain the terms with $\Phi=\Phi_0$. 
By symmetry, we will mainly restrict our attention to the region where 
\EQ{
 |\y|\ge|\y-\x|=m.} 
Since we are going to integrate by parts twice, we need up to the third derivatives of the phases. Let $I(r):=H''(r)/r-H'(r)/r^2$. Explicit computations give us 
\EQ{ \label{deris} 
 &H'(r) = \frac{2(1+r^2)}{\BR{r}},\quad 
  H''(r) = \frac{2r(3+r^2)}{\BR{r}^3},\quad 
  H'''(r) = \frac{12}{\BR{r}^5},\\  
 &H''''(r) = -\frac{60r}{\BR{r}^7},\quad 
  I(r) = -\frac{4}{r^2[r]^3},\quad 
  I'(r)=\frac{4(4+5r^2)}{r^3\BR{r}^5},}
As for the differences, we have for any $r\ge s\ge0$,
\EQ{ \label{diffs}
 &H(r)-H(s) \sim \LR{r}(r-s),\quad 
  H'(r)-H'(s) \sim \frac{r(r-s)}{\LR{r}},\\ 
 &|H''(r)-H''(s)| \lec \frac{r-s}{\LR{r}},\quad 
  |H'''(r)-H'''(s)| \lec \frac{r}{\LR{r}^2\LR{s}^5}(r-s),\\ 
 &\left|\frac{H'(r)}{r}-\frac{H'(s)}{s}\right| \lec \frac{r-s}{rs\LR{s}^3},\quad 
  |I(r)-I(s)| \lec \frac{r-s}{rs^2\LR{s}^3}.}
For any vector $v$, we denote the partial derivative with respect to $\y$ in the direction $v$ by $\p_v F(\x,\y) := v\cdot\p_\y F(\x,\y)$. 
We will omit the estimate with $\Phi_-$, which is easier than that with $\Phi_+$. For the phases $\Phi_0$ and $\Phi_+$, and for any vectors $a,b,c$, we have 
\EQ{
 \Phi &= H(\x)\pm H(\y)- H(\y-\x),\\
 \p_a\Phi &= \pm H'(\y)a_\y - H'(\y-\x)a_{\y-\x},\\
 \p_a\p_b\Phi &= \pm H''(\y)a_\y b_\y - H''(\y-\x)a_{\y-\x} b_{\y-\x} \\
 &\qquad \pm \frac{H'(\y)}{|\y|}a\pp\y\cdot b\pp\y - \frac{H'(\y-\x)}{|\y-\x|}a\pp{\y-\x}\cdot b\pp{\y-\x},\\
 \p_a\p_b\p_c\Phi &= \pm H'''(\y)a_\y b_\y c_\y - H'''(\y-\x)a_{\y-\x} b_{\y-\x}c_{\y-\x} \\
 &\qquad \pm I(\y) (a,b,c)_\y - I(\y-\x)(a,b,c)_{\y-\x},}
where the upper and lower signs correspond to $\Phi_0$ and $\Phi_+$ respectively, and $(a,b,c)_\y$ denotes the symmetric 3-tensor defined by 
\EQ{
 (a,b,c)_\y := a_\y b\pp\y \cdot c\pp\y + b_\y c\pp\y \cdot a\pp\y + c_\y a\pp\y \cdot b\pp\y.} 
We will use the following elementary geometry. For any $r\ge s\ge 0$ and unit vectors $\al,\be$, we have
\EQ{ \label{vec diff}
 |r\al-s\be| \sim |r-s| + s|\al-\be|,}
where $\le$ follows just by the triangle inequality and $\gec$ follows by squaring the both sides. For any nonzero vectors $a,b$, we have
\EQ{
 |\hat{a}+\hat{b}|^2 = 2(1+\hat{a}\cdot\hat{b}) = \frac{|a+b|^2-(|a|-|b|)^2}{|a||b|} = \frac{(|a|+|b|)^2-|a-b|^2}{|a||b|}.}
Hence by putting $(a,b)=(\y,\x-\y)$ and $(\y,\y-x)$, we have
\EQ{ \label{angl}
 &|\hat{\y}-\hat{\y-\x}|^2 = 2(1-\hat{\y}_{\y-\x}) = \frac{|\x|^2-(|\y|-|\y-\x|)^2}{|\y||\y-\x|},\\
 &|\hat{\y}+\hat{\y-\x}|^2 = 2(1+\hat{\y}_{\y-\x}) = \frac{(|\y|+|\y-\x|)^2-|\x|^2}{|\y||\y-\x|}.}

\subsection{Estimate for $|u|^2$, the case of $\Phi_0$}
First we consider the phase $\Phi_0$, for which the integration in $s$ does not play any role. For a fixed $0<\de\ll 1$ and each $\x\not=0$, we split the integral region of $\y$ into the following three overlapping domains:
\EQ{ \label{division Phi0}
 D_+(\x):=\{\y\in\R^d \mid &|\y|-|\y-\x|>(1-2\de)|\x|\},\\
 D_0(\x):=\{\y\in\R^d \mid 
 &||\y|-|\y-\x|| < (1-\de)|\x|\},\\
 D_-(\x):=\{\y\in\R^d \mid &|\y|-|\y-\x|<-(1-2\de)|\x|\}}
and choose a partition of unity $1=\chi_+(\y)+\chi_0(\y)+\chi_-(\y)$ satisfying $\supp\chi_*\subset D_*$ for $*=+,0,-$ and 
\EQ{ \label{cutoffbound1}
 |\na_\y^k\chi_*(\y)|\lec (|\y|+|\y-\x|)^{-k} \quad(|k|\le 2).}
Such functions can be given in the form $\chi((|\y|-|\y-\x|)/|\x|)$ 
with some one-dimensional cut-off function $\chi$, then its $\y$ derivatives are given by  
\EQ{
 &\p_a \chi = \frac{a_\y-a_{\y-\x}}{|\x|}\chi',\\
 &\p_a\p_b\chi = \frac{(a_\y-a_{\y-\x})(b_\y-b_{\y-\x})}{|\x|^2}\chi'' +
  \left[\frac{a^\perp_\y\cdot b^\perp_\y}{|\y|} - \frac{a^\perp_{\y-\x}\cdot b^\perp_{\y-\x}}{|\y-\x|}\right]\frac{\chi'}{|\x|}.}
The above bounds \eqref{cutoffbound1} follows from these identities, \eqref{angl}, and $|\x|\lec |\y|\sim|\y-\x|$ in $D_\pm\cap D_0$. 

By symmetry, it suffices to estimate only in $D_+$ and $D_0$. We first consider $D_+$. Here we use the polar coordinates $\y=r\th$, or in other words, we choose the direction $\hat{\y}=\th$ for the partial integration. 
By the definition of $D_+$, we have $|\y-\x|\le |\y|\gec|\x|$, and by \eqref{angl}, 
\EQ{ \label{angle est}
 1-\hat{\y}_{\y-\x} \lec \frac{\de |\x|^2}{|\y||\y-\x|}.}
Partial integration in $r$ gives 
\EQ{ \label{def Y}
 &\int F e^{i\Phi t} d\y = \frac{i}{t}\int (YF) e^{i\Phi t} d\y
 =\frac{-1}{t^2}\int (Y^2 F) e^{i\Phi t} d\y,}
where the operator $Y$ is defined by
\EQ{
 &Y:=\frac{1}{r}\p_r\frac{r}{\p_r\Phi}=\frac{1}{r\p_r\Phi}-\frac{\p_r^2\Phi}{(\p_r\Phi)^2} + \frac{\p_r}{\p_r\Phi},\\
 &Y^2=3\frac{(\p_r^2\Phi)^2}{(\p_r\Phi)^4}-3\frac{\p_r^2\Phi}{r(\p_r\Phi)^3}-\frac{\p_r^3\Phi}{(\p_r\Phi)^3}+\frac{2\p_r}{r(\p_r\Phi)^2}-3\frac{\p_r^2\Phi\p_r}{(\p_r\Phi)^3} + \frac{\p_r^2}{(\p_r\Phi)^2}.}
For the phase derivatives, we have the following estimates 
\EQ{
 \p_r\Phi_0 &= H'(\y)-H'(\y-\x) + H'(\y-\x)(1-\hat{\y}_{\y-\x})\\
 &\sim \frac{|\y|}{\LR{\y}}(|\y|-|\y-\x|) + \LR{\y-\x}(1-\hat{\y}_{\y-\x}) 
 \gec \frac{|\y||\x|}{\LR{\y}},\\
 \p_r^2\Phi_0 &= H''(\y)-H''(\y-\x) + \left[H''(\y-\x)-\frac{H'(\y-\x)}{|\y-\x|}\right]|\hat{\y}\pp{\y-\x}|^2\\
 |\p_r^2\Phi_0| &\lec \frac{|\y|-|\y-\x|}{\LR{\y}} + \frac{\LR{\y-\x}}{|\y-\x|}(1-\hat{\y}_{\y-\x}) \lec \frac{\p_r\Phi_0}{|\y-\x|},\\
 \p_r^3\Phi_0 &= H'''(\y)-H'''(\y-\x) + H'''(\y-\x)(1-\hat{\y}_{\y-\x}^3)\\ 
 &\qquad + 3I(\y-\x)\hat{\y}_{\y-\x}|\hat{\y}\pp{\y-\x}|^2\\
 |\p_r^3\Phi_0| &\lec \frac{|\y|-|\y-\x|}{\LR{\y-\x}^5\LR{\y}} + \frac{1-\hat{\y}_{\y-\x}}{|\y-\x|^2\LR{\y-\x}^3} \lec \frac{\p_r\Phi_0}{|\y-\x|^2},}
where we have used $|\y|\ge|\y-\x|$. Thus we obtain
\EQ{ \label{est D+}
 |Y^2(\chi_+ F)| &\lec \frac{|F|}{(\p_r\Phi_0)^2|\y-\x|^2} + \frac{|(\chi_+ F)_r|}{(\p_r\Phi_0)^2|\y-\x|} + \frac{|(\chi_+ F)_{rr}|}{(\p_r\Phi_0)^2}\\
 &\lec \frac{|\y||\y-\x|\LR{\y}^{2-\s}}{|\y|^2|\x|^2|\y-\x|^2\LR{\x}\LR{\y-\x}^{1+\s}}f(\x-\y)g(\y) \\
 &\sim \frac{\LR{\y}^{2-\s}}{|\x|^2\LR{\x}|\y||\y-\x|\LR{\y-\x}^{1+\s}}f(\x-\y)g(\y). }

Next we consider $D_0$. Here the main part of $\na\Phi_0$ is its angular component, which is not always close to either $\hat{\y}$ or $\hat{\y-\x}$, so we simply integrate in the direction of $\na\Phi_0$.\footnote{On the other hand, our method for estimating in $D_0$ is not adequate in $D_+$, where the angular difference terms such as $|H'(\y)/|\y|-H'(\y-\x)/|\y-\x||$ in \eqref{na bound} can not be controlled.} Partial integration gives 
\EQ{ \label{def A}
 &\int F e^{i\Phi t} d\y = \frac{i}{t} \int (AF)e^{i\Phi t} d\y = \frac{-1}{t^2}\int (A^2 F) e^{i\Phi t} d\y,}
where the operator $A$ is defined by 
\EQ{
 &A F := \na\cdot\frac{\na\Phi}{|\na\Phi|^2} F= \frac{\na\Phi}{|\na\Phi|^2}\cdot \na F - \frac{\na\Phi\cdot(2\na^2\Phi-\Delta\Phi I)\cdot \na\Phi}{|\na\Phi|^4} F,}
which satisfies
\EQ{ \label{na bound}
 &|A F| \lec \frac{|\na F|}{|\na\Phi|} + \frac{|\na^2\Phi|}{|\na\Phi|^2}|F|,\\
 &|A^2 F| \lec \frac{|\na^2 F|}{|\na\Phi|^2} + \frac{|\na^2\Phi|}{|\na\Phi|^3}|\na F| + \left[\frac{|\na^2\Phi|^2}{|\na\Phi|^4} + \frac{|\na^3\Phi|}{|\na \Phi|^3}\right]|F|.}
In the region $D_0$, we have $|\y|\sim|\y-\x|\gec|\x|$. Using \eqref{vec diff} together with \eqref{deris} and \eqref{diffs}, we have  
\EQ{
 |\na\Phi_0| &\sim |H'(\y)-H'(\y-\x)| + H'(\y-\x)|\hat{\y}-\hat{\y-\x}|\\ 
 &\gec \LR{\y}|\hat{\y}-\hat{\y-\x}|,\\
 |\na^2\Phi_0| &\lec |H''(\y)-H''(\y-\x)| + \left|\frac{H'(\y)}{|\y|}-\frac{H'(\y-\x)}{|\y-\x|}\right|\\
 &\quad+\left[H''(\y)+\frac{H'(\y)}{|\y|}\right]|\hat{\y}\otimes\hat{\y}-\hat{\y-\x}\otimes \hat{\y-\x}|\\
 &\lec \frac{|\x|}{|\y|^2} + \frac{\LR{\y}}{|\y|}|\hat{\y}-\hat{\y-\x}|,\\
 |\na^3\Phi_0| &\lec |H'''(\y)-H'''(\y-\x)| + |I(\y)-I(\y-\x)|\\
 &\quad + (H'''(\y)+|I(\y)|)|\hat{\y}-\hat{\y-\x}|\\
 &\lec \frac{|\x|}{|\y|^3\LR{\y}} + \frac{|\hat{\y}-\hat{\y-\x}|}{|\y|^2\LR{\y}^3}.}
By the definition of $D_0$ and \eqref{angl}, we have
\EQ{
 |\hat{\y}-\hat{\y-\x}|^2 = \frac{|\x|^2-(|\y|-|\y-\x|)^2}{|\y||\y-\x|} \sim \frac{|\x|^2}{|\y|^2}}
Hence we deduce that
\EQ{
 |\na\Phi_0|\gec \frac{\LR{\y}|\x|}{|\y|},\quad
 |\na^2\Phi_0| \lec \frac{|\na\Phi_0|}{|\y|},\quad
 |\na^3\Phi_0| \lec \frac{|\na\Phi_0|}{|\y|^2},}
and therefore
\EQ{
 |A^2(\chi_0F)| &\lec \frac{|\y||\x-\y|}{|\na\Phi_0|^2|\y|^2\LR{\x}\LR{\x-\y}^{1+\s}\LR{\y}^\s}f(\x-\y)g(\y)\\ 
 &\lec \frac{|\y|^2}{|\x|^2\LR{\x}\LR{\y}^{3+2\s}}f(\x-\y)g(\y),}
which is slightly better than the bound in \eqref{est D+}. 

In conclusion, we obtain
\EQ{ \label{t2 bound}
 t^2\left|\int_{|\y|\ge|\x-\y|} Fe^{i\Phi_0 t} d\y\right| &\lec \int_{1\gg|\y|\gec|\y-\x|} \frac{f(\x-\y)g(\y)}{|\x|^2|\y||\y-\x|} d\y\\
 &\qquad\qquad +\int_{1\lec|\y|\sim|\y-\x|} \frac{f(\x-\y)g(\y)}{|\x|^2\LR{\x}\LR{\y}^{1+2\s}} d\y\\
 &\qquad\qquad +\int_{1\lec|\y|\gg|\y-\x|} \frac{f(\x-\y)g(\y)}{\LR{\x}^{2+\s}|\y-\x|\LR{\y-\x}^{1+\s}} d\y\\
 &\lec \frac{\max(-\log|\x|,0)}{|\x|^2}\|f\|_{L^\infty}\|g\|_{L^\infty} + \frac{\|f\|_{L^2}\|g\|_{L^2}}{|\x|^2\LR{\x}^{2+2\s}}\\
 &\qquad + \LR{\x}^{-2-\s}\left[\frac{f(\x)}{|\x|\LR{\x}^{1+\s}}*g(\x)\right],}
where we used the Schwarz inequality and the condition $\s\ge -1/2$ for the second integral. For $1\gg|\x|$, we have also
\EQ{ \label{t0 bound}
 &\int |F|d\y \lec f'*g'.}
Applying this estimate in the region $\{|\x|<1/t\}$ and \eqref{t2 bound} in the rest, using the Young inequality $L^1*L^2\subset L^2$, $L^2*L^2\subset L^\infty$, and appending the same estimate in the opposite region $|\y-\x|\ge|\y|$, we obtain
\EQ{ \label{bd Phi0}
 &\norm{\LR{\x}^{1+\s}|\x|\int Fe^{i\Phi_0 t} d\y}_{L^2_\x}\\
 &\lec t^{-2}(\log t)^2\|f\|_{\FS}\|g\|_{\FS} + t^{-2}\|f'\|_{L^2}\|g'\|_{L^2},}
for $t\ge 2$. After integration in $t$, this estimate corresponds to the case $\th=1$ in \eqref{main est}. The remaining case $0<\th<1$ is covered by interpolation, see Section \ref{ss:interp}. 

\subsection{Estimate for $u^2$, the case of $\Phi_+$}
\begin{figure}[htb]
\hspace*{-16mm}
\epsfig{file=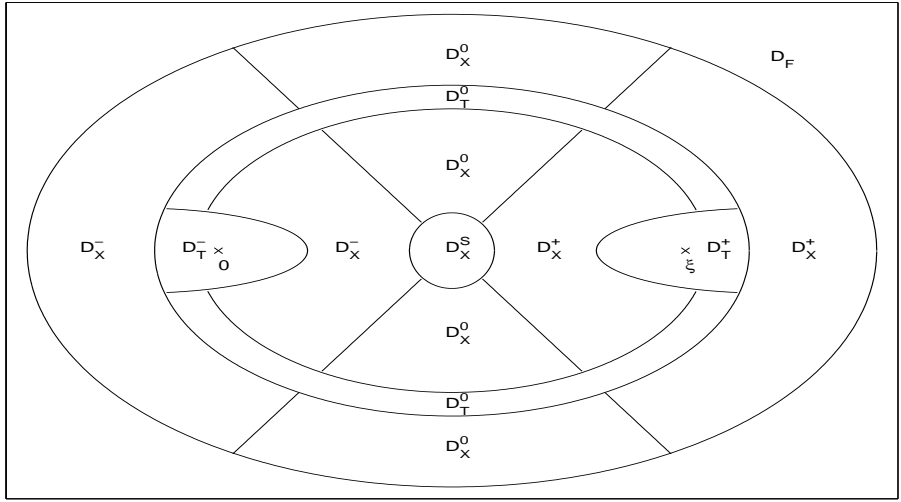, width=380pt, height=190pt}
\caption{Decomposition of $\eta$ space for $\Phi_+$ when $|\xi|\lec 1$ and
$t\gg 1$.}
\label{fig:1}
\end{figure}

Next we consider the phase $\Phi_+$, for which we need to take account of the time oscillation, and so split the integral into more regions. Let $\chi\in C^\infty(\R)$ satisfy $\chi(s)=1$ for $s\le 1$ and $\chi(s)=0$ for $s\ge 2$, and denote $\ti\chi=1-\chi$. Hence we have $\supp\chi\subset(-\infty,2]$ and $\supp\ti\chi\subset[1,\infty)$. We also denote 
\EQ{
  &\z=\y-\x/2,\quad r=|\y|,\quad \la=|\y|+|\y-\x|-|\x|,\\
  &M=\max(|\y|,|\y-\x|),\quad m=\min(|\y|,|\y-\x|).}
For a fixed positive $\de\ll 1$, and each $\x\not=0$ and $t\gg 1$, we introduce partitions of unity for $\y\in\R^d$ by the following identities: 
\EQ{
 1 = \chi_F + \chi_C,\quad \chi_C=\chi_T^++\chi_T^-+\chi_T^0 +\chi_X,\quad
 \chi_X=\chi_X^S+\chi_X^++\chi_X^-+\chi_X^0,}
and 
\EQ{
 &\chi_F = \ti\chi\left(\de \frac{\la}{|\x|}\right),\quad
 \chi_T = \chi_C\cdot\chi\left(\de \frac{\LR{\x}}{|\x|^{3}}\la\right), \\ 
 &\chi_T^\pm = \chi_T\cdot\chi_\pm, \quad \chi_T^0 = \chi_T\cdot\chi_0\cdot\ti\chi\left(\frac{\LR{\x}^2\la}{|\x|^{3}\de}\right),\\
 &\chi_X^S = \chi_X\cdot\chi(t|\z|),\ \chi_X^\pm = \chi_X\cdot\ti\chi(t|\z|)\cdot\ti\chi(\pm 4\hat{\z}_\x),}
where $\chi_\pm$ and $\chi_0$ are the same as in \eqref{division Phi0}. Hence denoting by $D_*^*=\supp\chi_*^*$, we have (see Figure \ref{fig:1})
\EQ{
 \R^d = D_F \cup D_T^+ \cup D_T^- \cup D_T^0 \cup D_X^S \cup D_X^+ \cup D_X^- \cup D_X^0.}
Remark that $D_X$ and $D_T$ cover the stationary-phase regions in space and time, respectively. $D_T$ is an annular region, separating $D_X$ into two connected components (if $|\x|\lec 1$). The derivatives of the cut-off functions satisfy 
\EQ{
 &|\na_\y^k\chi_F| \lec M^{-|k|},\quad 
 |(\hat{\y}\cdot\na_\y)^k\chi_T^*| + |(\hat{\y-\x}\cdot\na_\y)^k\chi_T^*| \lec m^{-k},\\
 &|\na_\y\chi_X^*| \lec |\z|^{-1} + D_T \LR{\x}|\x|^{-3/2}m^{-1/2},}
for $|k|\le 2$, where $D_T$ is identified with its characteristic function. We can easily derive these bounds using 
\EQ{
 &\p_r\la = 1+\hat{\y}_{\y-\x}\sim\frac{\la}{m},\quad 
 \p_r^2\la = \frac{|\hat{\y}^\perp_{\y-\x}|^2}{|\y-\x|}\lec\frac{\la}{m^2},\\
 &|\na_\y\la| = |\hat{\y}+\hat{\y-\x}| \sim \frac{\la^{1/2}}{m^{1/2}},\quad 
 |\na\chi_T^*| \lec \frac{\LR{\x}}{|\x|^{3/2}m^{1/2}},}
together with similar estimates for the radial derivative in $\y-\x$. 

In $D_F$, $D_T^+$ and $D_T^0$, we use the polar coordinates $\y=r\th$ and integrate twice in the radial direction.\footnote{Here we can not use the operator $A$ because there is no cancellation between the angular terms $H'(\y)/|\y|$ and $H'(\y-\x)/|\y-\x|$ in $\na^2\Phi_+$, and also because the cut-off function has larger derivative in the angular direction.} Since $|\y|\gec|\y-\x|+|\x|$, we have 
\EQ{
 -\p_r\Phi_+ &= H'(\y)-H'(\y-\x)+H'(\y-\x)(1+\hat{\y}_{\y-\x}),\\ 
 |\p_r^2\Phi_+| &\lec H''(\y)+H''(\y-\x)+\frac{H'(\y-\x)}{|\y-\x|}|\hat{\y}^\perp_{\y-\x}|^2\\
 &\lec \frac{|\y|}{\LR{\y}} + \frac{\LR{\y-\x}}{|\y-\x|}(1+\hat{\y}_{\y-\x}),\\ 
 |\p_r^3\Phi_+| &\lec H'''(\y)+H'''(\y-\x)+|I(\y-\x)||\hat{\y}^\perp_{\y-\x}|^2\\
 &\lec \frac{1}{\LR{\y-\x}^5} + \frac{1+\hat{\y}_{\y-\x}}{|\y-\x|^2\LR{\y-\x}^3}.}
In $D_F$, we have $|\y|\sim|\y-\x|\gg|\x|$, which implies that 
\EQ{
 -\p_r\Phi_+ \gec -\frac{|\y||\x|}{\LR{\x}} + \LR{\y} \gec \LR{\y}}
and also
\EQ{ \label{deri bounds} 
 |\p_r^{1+k}\Phi_+|\lec |\p_r\Phi|m^{-k}\quad (k=1,2).}
Therefore, defining $Y$ by \eqref{def Y} with $\Phi=\Phi_+$, we get 
\EQ{
 &|Y^2(\chi_{F}F)| \lec \frac{D_F |\y||\y-\x|f(\x-\y)g(\y)}{(\p_r\Phi_+)^2|\y-\x|^2\LR{\x}\LR{\y-\x}^{1+\s}\LR{\y}^\s} \lec \frac{f(\x-\y)g(\y)}{\LR{\x}^{4+2\s}}.} 
Since $\s>-1$, we obtain by the Schwarz inequality 
\EQ{ \label{est DF}
 \norm{\LR{\x}^{2+\max(\s,2\s)}\int \chi_F F e^{i\Phi_+t} d\y}_{L^2_\x}
 \lec t^{-2}\|f\|_{L^2}\|g\|_{L^2}.}
 
In $D_T^+$, we have $|\y|\sim|\x|\gec|\y-\x|$ and by the radial component 
\EQ{
 -\p_r\Phi_+ \gec \frac{|\y|^2}{\LR{\y}} \sim \frac{|\x|^2}{\LR{\x}},}
whereas in $D_T^0$ we have $|\y|\sim|\x|\sim|\y-\x|$ and by the angular component 
\EQ{
 -\p_r\Phi_+ \gec \LR{\x}(1+\hat{\y}_{\y-\x}) \gec \frac{|\x|^2}{\LR{\x}},}
and we have \eqref{deri bounds} in both cases. Hence we have for $*=0,+$, 
\EQ{
 |Y^2(\chi_{T}^*F)| &\lec \frac{D_T\LR{\x}^2|\y||\y-\x| f(\x-\y)g(\y)}{|\x|^4|\y-\x|^2\LR{\x}\LR{\y-\x}^{1+\s}\LR{\y}^\s}
 \lec \frac{D_T\LR{\x}^{1-\s}f(\x-\y)g(\y)}{|\x|^3|\y-\x| \LR{\y-\x}^{1+\s}}.}
For $|\x|\ge 1$, this gives the same bound as in the last term in \eqref{t2 bound}. For $|\x|\le 1$, we use the shape of $D_T$ and the polar coordinates for $\y-\x$, 
\EQ{
 \int|Y^2(\chi_{T}^*F)|d\y &\lec \int_0^{|\x|} \frac{\min(m,|\x|^3)\|f\|_{L^\infty}\|g\|_{L^\infty}}{|\x|^3 m} dm\\
 &\lec \|f\|_{L^\infty}\|g\|_{L^\infty}\log(1/|\x|),}
where the factor $\min(m,|\x|^3)$ is coming from integration in $\hat{\y-\x}$. 
We can treat $D_T^-$ in the same way by symmetry. Thus we obtain for $*=0,+,-$,  
\EQ{ \label{est DT}
 \norm{\LR{\x}^{2+\s}\int \chi_T^* F e^{i\Phi_+t} d\y}_{L^2_\x}
 \lec t^{-2}\|f\|_{L^2\cap L^\infty}\|g\|_{L^2\cap L^\infty}.}
 
Next we exploit the time oscillation in $D_X$, which is split into the following two cases:
\EQ{
 \begin{cases}
  \la\gg|\x|^3/\LR{\x},\\
  \la\ll|\x|^3/\LR{\x}^2,\quad||\y|-|\y-\x||\le (1-\de)|\x|.
 \end{cases}}
In the first case, we have $|\y|+|\y-\x|\sim|\x|\lec 1$ and 
\EQ{
 -\Phi_+ &= \sqrt{2}\la + \frac{|\y|^3}{[\y]+\sqrt{2}} + \frac{|\y-\x|^3}{[\y-\x]+\sqrt{2}} - \frac{|\x|^3}{[\x]+\sqrt{2}}
 \gec \frac{|\x|^3}{\LR{\x}},}
while in the second case, we have $|\y|\sim|\y-\x|\sim|\x|$ and 
\EQ{
 \Phi_+ &= [H(\x)-H(|\y|+|\y-\x|)] + [H(|\y|+|\y-\x|)-H(\y)-H(\y-\x)]\\
 &\gec \frac{|\x|^3}{\LR{\x}},}
where we used \eqref{diffs} and the identity 
\EQ{
 H(a+b)-H(a)-H(b) = \frac{ab(2a+b)}{\BR{a+b}+\BR{a}} + \frac{ab(a+2b)}{\BR{a+b}+\BR{b}} \quad(a,b>0).}
Therefore we can integrate in $s$ for $\x\not=0$:
\EQ{
 \X* := \left|\intS{*}\right| = \left|\int \frac{e^{i\Phi_+ t}}{i\Phi_+} (\chi_{X}^*F) d\y\right|,}
for $*=S,+,0,-$. 
In particular, we have
\EQ{
 \X{S} &\lec \int_{|\z|\lec 1/t} \frac{\LR{\x}}{|\x|^3} \frac{|\x|^2}{\LR{\x}^{2+2\s}} f(\x-\y)g(\y) d\y \\
 &\lec \frac{1}{|\x|\LR{\x}^{1+2\s}}\int_{|\z|\lec 1/t} f(\x/2+\z)g(\x/2-\z) d\z.}
Hence by using the triangle inequality we obtain 
\EQ{
 \|\LR{\x}^{2+2\s}\X{S}\|_{L^2_\x} &\lec \int_{|\z|\lec 1/t}\norm{\frac{\LR{\x}^{1/2}}{|\x|^{1/2}}f(\x/2+\z)}_{L^2_\x} \frac{\LR{\z}^{1/2}}{|\z|^{1/2}}\|g\|_{L^\infty} d\z\\
 &\lec t^{-3/2}\|f\|_{\FS}\|g\|_{L^\infty}.}
In $D_X^+$, we use the polar coordinate $\z=le^{i\om}\hat{\x}$ and partial integration in $l$: 
\EQ{
 \int F e^{i\Phi t} d\z = \frac{i}{t} \int (ZF) e^{i\Phi t} d\z,}
where the operator $Z$ is defined by
\EQ{
 Z := \frac{1}{l}\p_l\frac{l}{\p_l\Phi} = \frac{1}{l\p_l\Phi} - \frac{\p_l^2\Phi}{(\p_l\Phi)^2} + \frac{\p_l}{\p_l\Phi}.}
The first derivative is given by  
\EQ{
 \p_l\Phi_+ &= H'(\y)\hat{\y}_\z + H'(\y-\x)\hat{\y-\x}_\z \\
 &=(H'(\y)-H'(\y-\x))\hat{\y}_\z + H'(\y-\x)(\hat{\y}_\z+\hat{\y-\x}_\z),}
where the radial component can be estimated by
\EQ{
 (H'(\y)-H'(\y-\x))\hat{\y}_\z &\gec \frac{|\y|}{\LR{\y}}(|\y|-|\y-\x|)\frac{(\x+2\z)\cdot\hat{\z}}{|\y|}\\
 &\gec \frac{|\y|^2-|\y-\x|^2}{\LR{\x}}\frac{\x\cdot\hat{\z}}{|\y|} 
 \gec \frac{|\x||\z|\hat{\z}_\x^2}{\LR{\x}},}
and the angular component by 
\EQ{
 H'(\y-\x)(\hat{\y}_\z+\hat{\y-\x}_\z)
 &=H'(\y-\x) \frac{|\y|+|\y-\x|}{2|\z|}(1+\hat{\y}_{\y-\x})\\
 &\sim \LR{\y-\x}\frac{|\y|}{|\z|}|\hat{\y}+\hat{\y-\x}|^2.}
To rewrite it in terms of $\om$, let $\al$ and $\be$ be the angles $\in [0,\pi]$ such that
\EQ{
 \cos\al=\hat{\y-\x}\cdot\hat{-\x},\quad \cos\be=\hat{\y}\cdot\hat{\x}.}
Then by the sine theorem and $\om\sim\sin\om$ by $|\om|<\pi/2$, we have 
\EQ{ \label{lower angle}
 |\hat{\y}+\hat{\y-\x}| &\gec \sin\al+\sin\be \gec \left(\frac{1}{|\y-\x|}+\frac{1}{|\y|}\right)|\z||\sin\om| 
 \gec \frac{|\z|}{m}|\om|,}
and so 
\EQ{
 H'(\y-\x)(\hat{\y}_\z+\hat{\y-\x}_\z)
 \gec \frac{\LR{m}|\x||\z|}{m^2}\om^2
 \gec \frac{\LR{\x}|\z|\om^2}{m}.}
Thus we obtain 
\EQ{
 |\p_l\Phi_+| \gec \frac{|\x||\z|}{\LR{\x}} + \frac{\LR{\x}|\z|\om^2}{m}
 = \frac{\LR{\x}|\z|}{m}\left(\om^2+\frac{m|\x|}{\LR{\x}^2}\right).}
The second derivative is estimated by
\EQ{
 |\p^2_l \Phi_+| &\lec H''(\y) + H''(\y-\x) + \frac{H'(\y)}{|\y|}|\hat{\z}^\perp_\y|^2 + \frac{H'(\y-\x)}{|\y-\x|}|\hat{\z}^\perp_{\y-\x}|^2\\
 &\lec \frac{|\x|}{\LR{\x}} + \frac{\LR{m}}{m}\left|\frac{|\y|}{|\z|}\hat{\y}^\perp_{\y-\x}\right|^2 \lec \frac{|\x|}{\LR{\x}} + \frac{\LR{m}|\y|^2}{m|\z|^2}(1+\hat{\y}_{\y-\x}).}
Then by using $|\y|\sim \max(m,|\z|)$ and the above estimate, we obtain
\EQ{
 |\p^2_l \Phi_+| \lec |\p_l \Phi_+| [m^{-1}+|\z|^{-1}].}
Thus we obtain 
\EQ{ \label{est wX+}
 |Z(\chi_X^+F)| &\lec |\Phi_+| w_X^+ D_X f(\x-\y)g(\y), \quad \text{where}\\
 w_X^+ :&= \frac{\LR{\x}}{|\x|^3}\frac{m}{\LR{\x}|\z|\left(\om^2+\frac{m|\x|}{\LR{\x}^2}\right)}\frac{|\x| m}{\LR{\x}^{1+\s}\LR{m}^{1+\s}}\left[\frac{1}{m}+\frac{1}{|\z|}+\frac{D_T\LR{\x}}{m^{1/2}|\x|^{3/2}}\right]\\
 &= \frac{m^2}{\LR{\x}^{1+\s}|\x|^2\LR{m}^{1+\s}|\z|\left(\om^2+\frac{m|\x|}{\LR{\x}^2}\right)}\left[\frac{1}{m}+\frac{1}{|\z|}+\frac{D_T\LR{\x}}{m^{1/2}|\x|^{3/2}}\right].}
Using $|\y|\sim\max(m,|\z|)$, we have for $|\x|\ge 1$, 
\EQ{
 w_X^+ \lec \frac{1}{|\x|^{3+\s}\LR{m}^{1+\s}} + \frac{1}{|\x|^{2+2\s}|\z|^2},}
the first term is treated as the last term of \eqref{t2 bound}, and the second term by using the triangle and the Schwarz inequalities 
\EQ{
 \norm{\int_{1\ge |\z|\sgec 1/t}\frac{f(\x-\y)g(\y)d\z}{|\z|^2}}_{L^2_\x}
 &\lec \int_{1\ge |\z|\sgec 1/t}\frac{d\z}{|\z|^2} \|f(\x/2-\z)\|_{L^2_\x} \|g\|_{L^\infty},\\
 &\lec \|f\|_{L^2}\|g\|_{L^\infty}\log t,\\
 \norm{\int_{|\z|> 1}\frac{f(\x-\y)g(\y)d\z}{|\z|^2}}_{L^2_\x}
 &\lec \sqrt{\int_{|\z|>1} \frac{d\z}{|\z|^4}\int|f(\x-\y)g(\y)|^2 d\y d\x}\\
 &\lec \|f\|_{L^2}\|g\|_{L^2}.}
For $|\x|\le 1$, we have 
\EQ{
 w_X^+ \lec \frac{1}{|\z|^2(\om^2+|\x|^2)} + \frac{m^{1/2}D_T}{|\x|^{9/2}|\z|},} 
and so in the polar coordinates $\z=(l,\om)$, 
\EQ{
 \int_{D_X} w_X^+d\y
 &\lec \int_{1/t}^{|\x|}\int_0^{2\pi}\frac{l dl d\om}{l^2(\om^2+|\x|^2)}
  + \int_0^{|\x|}\frac{\min(l,|\x|^3)}{|\x|^4 l} dl\\
 &\lec\frac{\log t+|\log|\x||}{|\x|} 
}
Thus we obtain
\EQ{ \label{est DX+}
 \|\LR{\x}^{2+\min(\s,2\s)}\X+\|_{L^2_\x(|\x|\ge 1/t)} \lec t^{-1}(\log t)^2 \|f\|_{L^\infty\cap L^2}\|g\|_{L^\infty\cap L^2}.}
We have the same bound for $\X-$ by symmetry. 

In $D_X^0$, we use the operator $A$ defined by \eqref{def A} with $\Phi=\Phi_+$. Here we have $|\x|\sim|\y|\sim|\x-\y|$ and $|\hat{\z}_\x|\le 1/2$, which implies by the same argument as in \eqref{lower angle} that
\EQ{
 |\hat{\y}+\hat{\y-\x}| \gec \frac{|\z|}{|\x|}.}
Hence we have
\EQ{
 &|\na_\y\Phi_+| \gec H'(\y-\x)|\hat{\y}+\hat{\y-\x}| \gec \frac{\LR{\x}|\z|}{|\x|},\quad 
 |\na_\y^2\Phi_+| \lec \frac{\LR{\x}}{|\x|} \lec \frac{|\na_\y\Phi_+|}{|\z|},}
and so
\EQ{
 |A(\chi_X^0 F)| &\lec |\Phi_+| D_X^0 w_X^0 f(\x-\y)g(\y),\quad \text{where}\\
 w_X^0 :&= \frac{\LR{\x}}{|\x|^3} \frac{|\x|}{\LR{\x}|\z|}\frac{|\x|^2}{\LR{\x}^{2+2\s}}\left[\frac{1}{|\z|} + \frac{\LR{\x}D_T}{|\x|^2}\right]
 \lec \frac{1}{\LR{\x}^{2+2\s}|\z|}\left[\frac{1}{|\z|} + \frac{\LR{\x}D_T}{|\x|^2}\right].}
This bound is better than \eqref{est wX+}, and so in the same way we obtain
\EQ{ \label{est DX0} 
 \|\LR{\x}^{2+2\s}\X0\|_{L^2_\x(|\x|\ge 1/t)} \lec t^{-1}(\log t) \|f\|_{L^\infty\cap L^2}\|g\|_{L^\infty\cap L^2}.}

On the other hand, we have for $|\x|\le 1/t\lec 1$, 
\EQ{
 \iint_t^{2t} |\chi_X F| ds d\y &\lec \int_t^{2t} \int_{|\y|\lec|\x|} |\x|^2 \|f\|_{L^\infty} \|g\|_{L^\infty} d\y ds\\
 & \lec t^{-3} \|f\|_{L^\infty} \|g\|_{L^\infty} }
Putting the pieces together, we obtain
\EQ{
 \norm{\LR{\x}^{2+\min(\s,2\s)}\int_\infty^t \int e^{i\Phi_+ s} F d\y ds}_{L^2_\x}
 \lec t^{-1}(\log t)^2\|f\|_\FS \|g\|_\FS.}
Here the condition $\s\ge-1/2$ is inevitable for $\Phi_+$. The difference from the case of $\Phi_0$ is the inbalance between the $|\x|$ from $U$ and the $1/|\z|$ from the partial integration, which costs one regularity. On the other hand, $\Phi_+$ is better at $\x=0$ than $\Phi_0$, which is non-oscillatory at $\x=0$. 

\subsection{Interpolating estimates} \label{ss:interp}
Thus we have proved \eqref{main est} in the case $\th=1$. On the other hand, we have for any $k\in\Z$, 
\EQ{
 \norm{\int |F(\x,\y)| d\y}_{L^2(2^k<|\x|<2^{k+1})} 
 &\lec 2^k\norm{\int |F(\x,\y)|d\y}_{L^\infty_\x} \lec 2^k\|f'\|_{L^2}\|g'\|_{L^2}.}
which implies that for any $t>0$,
\EQ{
 \norm{\int_{2t}^{t} \int |\x|^{-1}|F(\x,\y)| d\y ds}_{L^2(2^k<|\x|<2^{k+1})}
 \lec t\|f'\|_{L^2}\|g'\|_{L^2}.}
This corresponds to the case $\th=-1$ in \eqref{main est}, although it is divergent for the integral $t\to\infty$. 
By applying real interpolation or the H\"older inequality to dyadic sequences on $|\x|$ and $t$, we obtain the desired estimate \eqref{main est} for $0<\th\le 1$. \qedsymbol
 
\section{Proof of the 2D theorem}
\subsection{Bootstrap setting}
We decompose our solution $z=z^0+z^1$ and $u=u^0+u^1$, where 
\EQ{
 z^0:=e^{-iHt}\fy,\quad u^0:=Vz^0,}
and starting from the above estimates, we will derive 
\EQ{ \label{goal}
 \|z^1\|_{\dot H^1} \lec t^{-\al},\quad \|z^1\|_{\dot H^{1/2}} \lec t^{-\be},\quad \|u\|_{L^4} \lec t^{-\be},}
with some $\al,\be$ satisfying
\EQ{ \label{cond ab}
 \be<1/2,\quad 1-\be<\al<2\be,}
by the standard iteration argument. $(\al,\be)$ can be arbitrarily close, but not equal, to $(1,1/2)$. We fix $\ka\in(0,1/4)$ such that
\EQ{
 1/2+\ka < \al < 2\be - \ka.}
Let $\|\fy\|_{\dot B^1_{1,1}}\le\de\ll 1$. The $L^p$ decay \eqref{decay est} implies the following bounds on the free part $u^0=u^0_1+iu^0_2$: 
\EQ{
 &\|u^0(t)\|_{H^1_x} \lec \|\fy\|_{H^1},\quad
  \|u_1^0\|_{L^\infty_{x}} + \|\na u^0\|_{L^\infty_{x}} \lec t^{-1}\|\fy\|_{\dot B^1_{1,1}} \le \de t^{-1},\\ 
 &\|u^0(t)\|_{L^4_x} \lec t^{-1/2}\|\fy\|_{\dot B^0_{4/3,2}}.}
The last quantity is finite for high frequency by interpolation of $\dot B^{1}_{1,1}$ and $\dot H^1$. The low frequency part is also finite, because \eqref{asy norm} implies that
\EQ{
 \LR{\na}^{-1/2}\fy \in\dot B_{p,\infty}^{2/p-2} \cap \dot B^{2/q-1}_{q,\infty},}
for all $p>1$ and $q\ge 1$, which is proved simply by estimating the inverse Fourier transform. 
In the following three subsections, we derive estimates on the normal form, the trilinear terms and the quadratic difference terms, where we need not assume that $u$ is the solution. For any function $u$, we denote  
\EQ{ 
 \|u\|_{\ZZ_T'} := &\|u\|_{\WL{\be}{0}{T}L^4_x}+\|\Re u\|_{\WL{\al}{0}{T}L^2_x}+\|\na u\|_{\WL{\al}{0}{T}L^2_x},\\ 
 \|u\|_{\ZZ_T} := &\|u^0\|_{\WL{1/2}{0}{T}L^4_x} +\|\Re u^0\|_{\WL{1}{0}{T}L^\infty_x}+\|\na u^0\|_{\WL{1}{0}{T}L^\infty_x} 
  +\|u-u^0\|_{\ZZ_T'}.}
Remark that $\ZZ_T$ is not a norm, but it is designed to measure different types of decay of $u^0$ and $u-u^0$, namely dispersive and dissipative. Since $\be<1/2$, we have
\EQ{
 \|u\|_{\WL{\be}{0}{T}L^4_x} \lec \|u\|_{\ZZ_T}.}

\subsection{Normal form} \label{ss:normal form}
The quadratic part is estimated just by the H\"older inequality: 
\EQ{ \label{est normal}
 &\||u|^2\|_{\WL{\al}{0}{T}L^2_x} \lec T^{-\ka}\|u\|_{\WL{\be}{0}{T}L^4_x}^2 \lec T^{-\ka}\|u\|_{\ZZ_T}^2,\\ 
 &\||u|^2-|w|^2\|_{\WL{\al}{0}{T}L^2_x}
  \lec T^{-\ka}(\|u\|_{\ZZ_T}+\|w\|_{\ZZ_T})\|u-w\|_{\ZZ_T'}.}

\subsection{Trilinear term} \label{ss:Tri}
For $1<p<\infty$, we have 
\EQ{
 \norm{N^3(u)}_{H^1_p}
 \lec \||u|^2u_1\|_{L^p} + \|u^2\na u\|_{L^p}.}
We apply this estimate after expanding $u=u^0+u^1$, choosing different $p$ for each term. For example, $(u_2)^2\na u_2$ is expanded into the following spaces
\EQ{
 &\|u_2 u_2 \na u_2^0\|_{\WL{2\be}{1}{T}L^2_x} \lec \|u_2\|_{\WL{\be}{0}{T}L^4_x}^2 \|\na u_2^0\|_{\WL{1}{0}{T}L^\infty_x},\\
 &\|u_2^0 u_2^0\na u_2^1\|_{\WL{\al+1/3}{1}{T} L^{2}_x} \lec \|u_2^0\|_{\WL{1/2}{0}{T}L^4_x}^{4/3} \|\na u_2^0\|_{\WL{1}{0}{T}L^\infty_x}^{2/3} \|\na u_2^1\|_{\WL{\al}{0}{T} L^2_x},\\
 &\|u_2 u_2^1 \na u_2^1\|_{\WL{2\al+\be-3/4}{3/4}{T}L^{4/3}_x} \lec \|u_2\|_{\WL{\be}{0}{T}L^4_x} \|\na u_2^1\|_{\WL{\al}{0}{T}L^2_x}^2,}
where we used the interpolation inequalities of Gagliardo-Nirenberg type:
\EQ{
 &\|u\|_{L^\infty_x} \lec \|u\|_{L^4_x}^{2/3} \|\na u\|_{L^\infty_x}^{1/3},\quad 
  \|u\|_{L^8_x} \lec \|u\|_{L^4_x}^{1/2}\|\na u\|_{L^2_x}^{1/2}.}
The other terms containing $u_1$ are estimated in the same way. Thus we obtain 
\EQ{ \label{est tri} 
 &\|\Tri(u)\|_{Stz_T^1} \lec T^{-\al-\ka}\|u\|_{\ZZ_T}^3,\\
 &\|\Tri(u)-\Tri(w)\|_{Stz_T^1} \lec T^{-\al-\ka}\left[\|u\|_{\ZZ_T}+\|w\|_{\ZZ_T}\right]^2\|u-w\|_{\ZZ_T'}. }

\subsection{Quadratic error term} \label{ss:Dif}
The quadratic difference term can be expanded by putting $u=u^0+u^1$ 
\EQ{
 &N^2(u^0+u^1)-N^2(u^0)\\
 &=-2i(2u^0_1+u^1_1)u^1_1-2PU^{-1}\na\cdot(u^0_1\na u^1_2 + u^1_1\na u^0_2 + u^1_1\na u^1_2),}
and each term is estimated in $T^{-\al}Stz_T^1$ by using
\EQ{
 &\|u^0_1u^1_1\|_{H^1} \lec \|u^0_1\|_{W^{1,\infty}}\|u^1_1\|_{H^1} \lec \de t^{-1-\al}\|u\|_{\ZZ_t},\\
 &\|u^1_1 u_1^1\|_{H^1_{p}} \lec \|u^1_1\|_{L^q}\|u^1_1\|_{H^1} \lec t^{-2\al}\|u\|_{\ZZ_t}^2,\\
 &\|u^0_1\na u^1_2\|_{L^2} \lec \|u^0_1\|_{L^\infty}\|\na u^1_2\|_{L^2} \lec \de t^{-1-\al}\|u\|_{\ZZ_t},\\
 &\|u^1_1\na u^0_2\|_{L^2} \lec \|u^1_1\|_{L^2}\|\na u^0_2\|_{L^\infty} \lec \de t^{-1-\al}\|u\|_{\ZZ_t},\\
 &\|u^1_1\na u^1_2\|_{L^{p}} \lec \|u^1_1\|_{L^q}\|\na u^1_2\|_{L^2} \lec t^{-2\al}\|u\|_{\ZZ_t}^2,}
where $\de$ is the small factor coming from $\|\fy\|_{\dot B^1_{1,1}}$, and we choose $p\in(1,2)$ and $q\in(2,\infty)$ such that
\EQ{
 \al-\ka>1-1/q=3/2-1/p.}
Thus we obtain 
\EQ{ \label{est dif}
 &\|\Dif(u)\|_{Stz_T^1} \lec T^{-\al}\left[\de+T^{-\ka}\|u\|_{\ZZ_T}\right] \|u\|_{\ZZ_T},\\ 
 &\|\Dif(u)-\Dif(w)\|_{Stz_T^1}
  \lec T^{-\al} \left[\de+T^{-\ka}(\|u\|_{\ZZ_T}+\|w\|_{\ZZ_T})\right] \|u-w\|_{\ZZ_T'}.}

\subsection{Iteration argument}
We define an iteration sequence $(z_{(k)},u_{(k)},v_{(k)})$ for $k=0,1,2\dots$ and $t>T\gg 1$, by
\EQ{ \label{itera}
 &z_{(0)}=v_{(0)}=z^0=e^{-iHt}\fy,\quad u_{(0)}=u^0=Vz^0,\\
 &z_{(k+1)} = z^0 + \Tri(u_{(k)}) + \Dif(u_{(k)}) + \Asy(u^0),\\
 &u_{(k+1)} = Vv_{(k+1)} = Vz_{(k)} - P\frac{|u_{(k)}|^2}{2}. }
We introduce the following norm for $z$:
\EQ{
  \|z\|_{\ZZ_T^2} := \|z\|_{\WL{\al}{0}{T}\dot H^1_x} + \|z\|_{\WL{\be}{0}{T}\dot H^{1/2}_x}.}
Since $\be<\al$, we have $\|z\|_{\ZZ_T^2}\lec T^\al\|z\|_{Stz_T^1}$. Using the Sobolev embedding $\dot H^{1/2}_x\subset L^4_x$ and $\be<\al$, we have for any $k,j=0,1,2\dots$, 
\EQ{
 &\|u_{(k+1)}-u_{(j+1)}\|_{\ZZ_T'} \sim \|v_{(k+1)}-v_{(j+1)}\|_{\WL{\al}{0}{T}\dot H^1_x} + \|u_{(k+1)}-u_{(j+1)}\|_{\WL{\be}{0}{T}L^4_x} \\ 
 &\lec \|z_{(k)}-z_{(j)}\|_{\ZZ_T^2} + \||u_{(k)}|^2-|u_{(j)}|^2\|_{\WL{\al}{0}{T}L^2_x},}
and then the quadratic part is estimated by using \eqref{est normal} and $\al+\ka<2\be$,
\EQ{
 \||u_{(k)}|^2-|u_{(j)}|^2\|_{\WL{\al}{0}{T}L^2_x}
 \lec T^{-\ka}\left[\|u_{(k)}\|_{\ZZ_T}+\|u_{(j)}\|_{\ZZ_T}\right]\|u_{(k)}-u_{(j)}\|_{\ZZ_T'}.}
For the first iteration, we have $u_{(1)}-u_{(0)}=-P|u^0|^2/2\in\R$, and so
\EQ{ \label{first u}
 \|u_{(1)}-u_{(0)}\|_{\WL{\be}{0}{T}L^4_x} \lec \|v_{(1)}-v_{(0)}\|_{\WL{\al}{0}{T}\dot H^1_x} \lec \||u^0|^2\|_{\WL{\al}{0}{T}L^2_x} \lec T^{-\ka}\|u^0\|_{\ZZ_T}^2.}
As for $z_{(k)}$, we use \eqref{est tri} and \eqref{est dif}, deriving 
\EQ{
 \|z_{(k+1)}-z_{(j+1)}\|_{Stz_T^1} \lec T^{-\al} &\left[\de+T^{-\ka}(1+\|u_{(k)}\|_{\ZZ_T}+\|u_{(j)}\|_{\ZZ_T})^2\right]\\
  &\times \|u_{(k)}-u_{(j)}\|_{\ZZ_T'}.}
For the first iteration, we apply \eqref{est tri} and \eqref{main est}. Then we get 
\EQ{
 \|z_{(1)}-z_{(0)}\|_{\ZZ_T^2}
  &\lec \|\Tri(u^0)\|_{T^{-\al}Stz_T^1}+\|\Asy(u^0)\|_{\ZZ_T^2}\\ 
  &\lec T^{-\ka}\left[\|u^0\|_{\ZZ_T}^3 + \|\fy\|_{\NN}^2\right],}
where we denote
\EQ{
 \|\fy\|_{\NN} := \|\fy\|_{H^1} + \sum_{0\le|k|\le 2} \|\LR{\x}^{-1/2}|\x|^{|k|}\p_\x^k\ti\fy(\x)\|_{L^2\cap L^\infty}.}
Gathering the above estimates, we deduce that
\EQ{
 &D_k :=\|z_{(k)}-z_{(k-1)}\|_{\ZZ_T^2} + \|u_{(k+1)}-u_{(k)}\|_{\ZZ_T'},\quad 
  E_k :=\|u_{(k)}\|_{\ZZ_T},\\ 
 &D_{k+1} \lec (\de+T^{-\ka}(1+E_{k+1}+E_k+E_{k-1})^2)D_{k},\quad(k=1,2,3\dots),\\
 &D_1 \lec T^{-\ka}\left[\|u^0\|_{\ZZ_T}^2(1+\|u^0\|_{\ZZ_T})+\|\fy\|_{\NN}^2\right],\\
 &E_k \lec \sum_{j=1}^{k-1} D_j + \|u_{(1)}\|_{\ZZ_T},\quad 
  \|u_{(1)}\|_{\ZZ_T} \lec \|u^0\|_{\ZZ_T}(1+T^{-\ka}\|u^0\|_{\ZZ_T}).}
Hence for sufficiently small $\de$ and large $T$, $(z_{(k)},u_{(k)})$ converges to some function $(z,u)$ satisfying the equation \eqref{eq nu} and 
\EQ{
 \|z-z_{(0)}\|_{\ZZ_T^2} + \|u-u_{(1)}\|_{\ZZ_T'} 
 + \|z-z_{(1)}\|_{T^{-\al}Stz_T^1} \lec D_1.}
The uniqueness for $t>T$ is proved also by the above difference estimates. In addition, \eqref{main est} implies that 
\EQ{
 &t^{1-\eps}(\|z'\|_{\dot H^1}+\|z_{(1)}-z^0-z'\|_{L^2_x}) 
   + t^{\eps/2}\|z'\|_{\dot H^\eps} 
 \lec \|\fy\|_{\NN}^2,}
where we used the fact that $z_{(1)}-z^0-z'=\Asy'(u^0)$ does not contain $\Phi_0$. By using this estimate together with \eqref{first u} and the $L^p$ decay, we also obtain
\EQ{
 \|\nu\|_{\dot H^2\cap\dot H^1} &\lec \|u\|_{L^4_x}^2 \lec t^{-2\be}(D_1+\|u^0\|_{\ZZ_T})^2,\\
 \|\nu\|_{\dot H^{2\eps}} &\lec \|u\|_{L^{2/(1-\eps)}_x}^2 
  \lec (\|z\|_{L^{2/(1-\eps)}_x} + \|u\|_{L^4_x}^2)^2\\
 &\lec t^{-\eps}(\|\fy\|_{\dot B^0_{2/(1+\eps),1}}+\|\fy\|_{\NN}^2+D_1^2+\|u^0\|_{\ZZ_T}^2)^2,}
for any small $\eps>0$. In particular, $z(T),v(T)\in\dot H^s$ for $0<s\le 1$. The local uniqueness for \eqref{eq nu T} in this class is easily derived from the Strichartz, Sobolev and H\"older inequalities. 

\subsection{Global continuation} \label{ss:GC}
The final task is to extend our solution $u$ to $t<T$. We can not apply the $H^1$ global wellposedness of $u$ by \cite{BS}, since $\Asy(u^0)$ barely falls out of $L^2$. However the nonlinear energy is still finite, because $u_1,|u|^2\in L^2$, and the $L^2$ singular part at low frequency belongs to $L^\infty$. Hence we can apply the global existence results in \cite{Ge}, and we have only to see persistence of our function space, namely $(z,v)\in C(\R;\dot H^\eps\cap\dot H^1)$. The conserved energy can be written as 
\EQ{
 E(u) := \int \frac{|\na\psi|^2}{2} + \frac{(|\psi|^2-1)^2}{4} dx = \int \frac{|\na u|^2}{2} + |\rh|^2 dx,}
where we denote $\rh:=u_1+\frac{|u|^2}{2}$. \cite[Theorem 1.1]{Ge} gives 
\EQ{
 &\rh\in C(\R;L^2),\quad u\in C(\R;(L^\infty + H^1)\cap\dot H^1),\\
 &u-e^{i\Delta(t-T)}u(T)\in C(\R;H^1).}
Hence we have $u \in C(\R;\dot H^s)$ for $0<s\le 1$, which implies for $0<s<1$ that  
\EQ{ \label{quad space} 
 PU^{-1}|u|^2 \in C(\R;\dot H^{2s}).} 
On the other hand, by using the identities 
\EQ{
 Uz_1 = u_1 + P\frac{|u|^2}{2} = P\rh + Q u_1,\quad z_2 = u_2,}
we get
\EQ{
 z\in C(\R;(L^\infty+H^1)\cap\dot H^1),\quad \na\rh\in C(\R;L^2+L^{8/5}),}
and the integral equation for $z$ can be written as 
\EQ{
 z = &e^{-iH(t-T)}z(T) \\
  &- \int_T^t e^{-iH(t-s)}[2i(\rh u_1) + 2PU^{-1}\na\cdot(u_1\na u_2) - 2U(\rh u_2-u_1u_2)] ds,}
where the nonlinearity is in $C(\R;H^1+H^{1,4/3})$. Hence the Strichartz estimate implies that 
\EQ{
 z(t)-e^{-iH(t-T)}z(T) \in C(\R;H^1),}
and therefore 
\EQ{
 z \in C(\R;\dot H^{s}),}
for $0<s\le 1$. Combined with \eqref{quad space}, this implies that
\EQ{
 v=z-PU^{-1}|u|^2/2 \in C(\R;\dot H^s),}
for $0<s\le 1$. 
\qedsymbol

\section*{Acknowledgments}
The research of Gustafson and Tsai is partly supported by NSERC grants. The research of Nakanishi was partly supported by the JSPS grant no.~15740086.

\bigskip

\noindent{Stephen Gustafson},  gustaf@math.ubc.ca \\
Department of Mathematics, University of British Columbia, 
Vancouver, BC V6T 1Z2, Canada

\bigskip

\noindent{Kenji Nakanishi},
n-kenji@math.kyoto-u.ac.jp \\
Department of Mathematics, Kyoto University,
Kyoto 606-8502, Japan

\bigskip

\noindent{Tai-Peng Tsai},  ttsai@math.ubc.ca \\
Department of Mathematics, University of British Columbia, 
Vancouver, BC V6T 1Z2, Canada

\end{document}